\documentclass[11pt,a4paper, twoside]{article}
\usepackage[notref,notcite]{showkeys}

\usepackage{amsmath,amssymb,amsthm,amsfonts,mathrsfs,amscd,environ}
\usepackage{dsfont}
\usepackage{latexsym,enumerate,color,geometry,extarrows}
\usepackage{xcolor}
\usepackage{verbatim,fancyhdr,enumitem}
\def\beg{\begin}
\def\bequ{\begin{equation}}
\def\enqu{\end{equation}}
\def\bes{\begin{split}}
\def\ens{\end{split}}
\def\bews{\begin{ews}}
\def\beqn{\begin{eqnarray}}
\def\enqn{\end{eqnarray}}
\def\beq*{\begin{equation*}}
\def\enq*{\end{equation*}}
\def\bqn*{\begin{eqnarray*}}
\def\eqn*{\end{eqnarray*}}
\def\bary{\begin{array}}
\def\eary{\end{array}}
\def\bpma{\begin{pmatrix}}
\def\epma{\end{pmatrix}}
\def\bvma{\begin{Vmatrix}}
\def\evma{\end{Vmatrix}}

 \numberwithin{equation}{section}

\def\al{\alpha}
\def\be{\beta}
\def\ga{\gamma}
\def\de{\delta}

\def\th{\theta}

\def\ka{\kappa}
\def\la{\lambda}

\def\si{\sigma}

\def\Ph{\Phi}

\def\Om{\Omega}

\def\Q{\mathbb Q}
\def\R{\mathbb R}
\def\P{\mathbb P}
\def\E{\mathbb E}

\def\sF{\mathscr F}

\def\sC{\mathscr C}
\def\sB{\mathscr B}

\def\sL{\mathscr L}

\def\sP{\mathscr P}

\def\cC{\mathcal C}

\def\cT{\mathcal T}

\def\d{\mathrm{d}}

\def\ff{\frac}
\def\ra{\rightarrow}
\def\nn{\nabla}
\def\pp{\partial}
\def\<{\langle}
\def\>{\rangle}
\def\sq{\sqrt}
\def\tld{\tilde}
\def\we{\wedge}
\def\1{\mathds{1}}

\def\var{\text {\rm Var}}

\allowdisplaybreaks

\title{{\bf Exponential convergence in Wasserstein metric for distribution dependent SDEs}
}

\author{
{\bf Shao-Qin Zhang }\\
\footnotesize{School of Statistics and Mathematics, Central University of Finance and Economics, Beijing 100081, China}\\
\footnotesize{Email: zhangsq@cufe.edu.cn}\\
}

\begin{document}

\maketitle

\begin{abstract}

\end{abstract}\noindent
The existence and uniqueness of stationary distributions and the exponential convergence in $L^p$-Wasserstein distance are derived for distribution dependent SDEs   from associated  decoupled equations.  To establish the exponential convergence, we introduce a twinned Talagrand inequality of the original SDE and the associated decoupled equation, and  explicit convergence rate is obtained. Our results can be applied to  SDEs without uniformly dissipative drift and distribution dependent diffusion term, which cover the Curie-Weiss model  and the granular media model in double-well landscape with quadratic interaction as examples.

\noindent


\vskip 2cm

\section{Introduction}

Stochastic differential equations with distribution dependent drifts were introduced  by McKean \cite{McK} to investigate Vlasov-Poisson-Fokker-Planck systems. These type SDEs have attracted great attention since then, see e.g.  \cite{Daw,Mel,Szn} and recent works \cite{BLPR,CarDe,CGPS,HRW,LWZ,Wan18,Wan21a} with references therein.  Let $\sP $ denote the space of probability measures on $\R^d$ equipped  with the weak  topology.  Consider the following distribution dependent SDE on $\R^d$
\beg{align}\label{equ_DD}
\d X_t=b(X_t,\sL_{X_t})\d t+\si(X_t,\sL_{X_t})\d B_t,~t\geq 0,
\end{align}
where $\{B_t\}_{t\geq 0}$ is a $d$-dimensional Brownian motion on a complete filtration probability space $(\Om,\sF,\{\sF_t\}_{t\ge 0},\P)$, $\sL_{X_t}^{\P}$ is the law of  $X_t$ under $\P$, and
\beg{align*}
b:~\R^d\times\sP\ra \R^d,\qquad \si:~\R^d\times \sP\ra \R^d\otimes\R^d
\end{align*}
are measurable. If $\si(x,\mu)=\si(x)$ is independent of $\mu$, \eqref{equ_DD} is also called the McKean-Vlasov SDE. If $b(x,\mu)=b(x)$ moreover, \eqref{equ_DD} becomes the classical time homogenous It\^o-type SDE. Distribution dependent SDEs can be derived from the associated interacting particles system by passing to the mean field limit, and the distribution dependent part of coefficients reflects the interaction of the particles system, see \cite{Szn} for example. The well-posedness  for \eqref{equ_DD} in the weak and strong sense has been intensively investigated, see e.g. \cite{HRW,RZ,Wan18,Wan21a,ZXC} and references within.   

The convergence to the equilibrium of the solution to McKean-Vlasov SDEs has been widely studied. In the case that $\si=\sq 2 I$ with $I$ the identity matrix  on $\R^d$ and $b(x,\mu)=-\nn V(x)-\nn F*\mu(x)$ where $V,F\in C^2(\R^d)$ with $F(-x)=F(x)$, $\nn $ is the gradient operator and  $*$ stands for the convolution on $\R^d$:
$$f*\mu(x)=\int_{\R^d} f(x-y)\mu(\d y),~f\in\sB(\R^d),$$ 
\cite{CMV} obtained the explicit exponential convergence in mean field entropy for \eqref{equ_DD} in a variety of convexity conditions on confining potential $V$ and  interaction potential $F$.  Recently, \cite{GLW} generalized results in \cite{CMV} by using functional inequalities and establishing detail estimates on the associated interaction particles system. \cite{Wan18} obtained  existence and uniqueness of stationary probability measures and the exponential convergence in Wasserstein distance for \eqref{equ_DD} with dissipative drifts and distribution dependent $\si$ satisfying 
\beg{align}\label{Wdis}
2\<b(x,\mu)-b(y,\nu),x-y\>&+\|\si(x,\mu)-\si(y,\nu)\|_{HS}^2\nonumber\\
&\leq C_1W_2(\mu,\nu)^2-C_2|x-y|^2,~x,y\in\R^d,\mu,\nu\in\sP^2.
\end{align}
By using the log-Harnack inequality and the Talagrand inequality, \cite{RenW} established the exponential convergence in classical entropy and Wasserstein distance under general setting of $b$ but  distribution-free $\si$, which extended researches of \cite{CGPS,GLW,Wan18}. For the general non-convex case, drift term $b$ is not uniformly dissipative w.r.t. the first variable, i.e.   \eqref{Wdis}  holds only for large $|x-y|$.  A quantitative method that combines Lyapunov functions with reflection coupling and concave distance functions is developed  to investigate the longtime behavior of McKean-Vlasov SDEs without uniformly dissipative drifts, see e.g. \cite{EGZ,LWZ,Wan21}.  In this paper, we consider \eqref{equ_DD} with  general distribution dependent $\si$ and  without uniformly dissipative drifts.  Many methods mentioned previously fail in this case.  

In our previous paper \cite{ZSQ21},  existence results on stationary probability measures  and criteria on phase transition (the existence of multi-stationary states)  have been investigated  for \eqref{equ_DD}. The phase transition can occur for the general non-convex case with strong interaction in particular, see \cite{ZSQ21} or  \cite{CGPS,Daw,FenZ,Tug14a} as well as references within. In this paper, the existence and uniqueness of stationary probability measures and the exponential convergence are established under estimates of the weakness of the  interaction, see $\de_0,\de_1$ and $\de_2$ in Theorem \ref{thm-uni}, Theorem \ref{thm-exp0} and Theorem \ref{thm-exp} below.


Given $\mu\in\sP$, there is a decoupled equation associated with \eqref{equ_DD}    
\beg{align}\label{equ_dec}
\d X_t^{\mu}= b(X_t^\mu,\mu)\d t+\si(X_t^\mu,\mu)\d B_t.
\end{align}
Since $\mu$ is fixed, \eqref{equ_dec} is a classical time homogenous It\^o-type SDE. When \eqref{equ_dec} is well-posed and $X_0^{\mu}=x\in\R^d$, we denote by $P^\mu_tf(x)=\E f(X_t^\mu)$, where $f\in\sB_b(\R^d)$ and $\sB_b(\R^d)$ consists of all bounded Borel measurable functions on $\R^d$. $P_t^\mu$ is the Markov semigroup associated with \eqref{equ_dec}. There are rich researches on classical It\^o-type SDEs to establish the exponential convergence for the solution to \eqref{equ_dec}, see e.g. \cite{Eb,EGZ,LuoW,W20} and references therein. We derive  existence and uniqueness of stationary distributions for \eqref{equ_DD} from the ergodicity of \eqref{equ_dec}.  For any $p\geq 1$, let
\beg{align*}
\sP^p&=\{\mu\in\sP ~|~\|\mu\|_p:=(\mu(|\cdot|^p))^{\ff 1 p}<\infty\}.
\end{align*}
We denote by $W_p$ the $L^p$-Wasserstein distance on $\sP^p$:
\beg{align*}
W_p(\mu,\nu):= \inf_{\pi\in\sC(\mu,\nu)}\left(\int_{\R^d\times\R^d}|x-y|^p\pi(\d x,\d y)\right)^{\ff 1 p},~\mu,\nu\in\sP^p,
\end{align*} 
where $\sC(\mu,\nu)$ consists of all couplings of $\mu$ and $\nu$. $\sP^p$ becomes a complete metric space under the distance $W_p$. We assume that
\beg{description}[align=left, noitemsep]
\item[(H)] There is $p\geq 1$ so that for every $\mu\in\sP^p$, $P^\mu_t$ has a unique invariant probability measure $\cT_\mu\in \sP^p$, and there are $\hat C>0,\hat\la>0$ independent of $\mu$ such that
\beg{align}\label{exp-c1}
W_p((P_t^\mu)^*\nu,\cT_\mu)\leq \hat C e^{-\hat \la t} W_p(\nu,\cT_\mu),~\nu\in\sP^p. 
\end{align}
\end{description}
Combining this with the Talagrand inequality for \eqref{equ_dec}, we prove the mapping $\mu\ra \cT_\mu$ is contractive on $\sP^p$ when $b,\si$ depend weakly on the distribution, see Theorem \ref{thm-uni} below.   For concrete conditions that ensure {\bf (H)}, one can see Remark \ref{rem0} or {\bf (A2)} with Corollary \ref{cor2} or {\bf (A2')} with Corollary \ref{cor3}. To establish exponential convergence,  besides using the Talagrand inequality of the stationary distribution (Theorem \ref{thm-exp0}), we introduce a twinned Talagrand inequality of \eqref{equ_DD} and \eqref{equ_dec}, see {\bf (Ta)} and Theorem \ref{thm-exp} for details.

This paper is structured as follows. Our main results and corollaries are stated in Section 2.  Proofs of main results are given in Section 3, and proofs of corollaries are given in Section 4.

\section{Main results and corollaries}
\subsection{Existence and uniqueness}
We first present the existence and uniqueness of  stationary distributions to \eqref{equ_DD}, which is also equivalent to that there is a unique $\mu$ such that \eqref{equ_dec} with initial distribution $\mu$ has a unique weak solution. Hence, we only need to concern with the weak well-posedness of \eqref{equ_dec}. We assume that $b,\si$ satisfy the following assumption, which also implies the strong well-posedness of \eqref{equ_dec} indeed.   
\beg{description}[align=left, noitemsep]
\item[(A1)] $b$ is continuous  in the  first variable, $\si$ is bounded and Lipschitz in the first variable, and there exist $K_0\in\R$ and $\de\geq 0$ such that  
\beg{align}\label{ine-bsi}
2\<b(x,\mu)-b(y,\nu),x-y\>&+(1+(p-2)^+)\|\si(x,\mu)-\si(y,\nu)\|_{HS}^2\nonumber\\
&\leq K_0|x-y|^2+\de^2W_p(\mu,\nu)^2,~\mu,\nu\in\sP^p.
\end{align}
\end{description}
To see the strong well-posedness of \eqref{equ_dec} with fixed $\mu$, we can set $\nu=\mu$ and $y=0$ in \eqref{ine-bsi}. Then, taking into account that $\si(x,\mu)$ is bounded in $x$,  there is $C>0$  such that 
\beg{align*}
2\<b(x,\mu),x\>+\|\si(x,\mu)\|_{HS}^2&\leq K_0|x|^2+2\<b(0,\mu),x\>+ \|\si(0,\mu)\|_{HS}^2\\
&\quad+\|\si(x,\mu)\|_{HS}\|\si(0,\mu)\|_{HS} \\
&\leq C(1+|x|^2).
\end{align*}
Combining this  with \eqref{ine-bsi} (setting $\nu=\mu$), it follows from Krylov's criterion, see  e.g. \cite[Theorem 3.1.1]{LiuR}, that \eqref{equ_dec} has a unique solution.

Before our first theorem, we give a simple proposition which indicates that under {\bf (H)}, \eqref{equ_DD} has a unique stationary probability measure for small $\de$.
\beg{prp}\label{prp1}
Assume that {\bf(H)} and ${\bf(A1)}$ hold. If $\de<\de_0$ with
$$\de_0:=\sup_{t>\hat\la^{-1}\log\hat C}\left\{\left(\ff {2t (1-\exp\{-(\ff {p\vee 2} 2 K_0+\ff {(p-2)^+} 2)t\})} {(p\vee 2)K_0+(p-2)^+}\right)^{-\ff 1 {p\vee 2}}   (1-\hat C e^{-\hat \la t})\right\},$$
then there is a  unique stationary probability measure for \eqref{equ_DD}.
\end{prp}

Let $H(\nu|\mu)$ be the relative entropy of $\nu$ with respect to $\mu$: 
\beg{align*}
H(\nu|\mu)=\begin{cases}\int_{\R^d} \log \ff{\d\nu}{\d\mu}\d\nu,\quad \ \ & 
 \mbox{if}~\nu\ll\mu,\\ 
  +\infty, \quad \ \ & \mbox{otherwise}.
\end{cases} 
\end{align*}
If  the invariant probability measure of $P_t^\mu$ satisfies the Talagrand inequality, then we have the following theorem. Let 
$$K(m,p)=2^{-\ff 1 {p\vee 2}}\left((p\vee 2)(2m-K_0)-(p-2)^+\right)^{\ff 1 {p\vee 2}}.$$

\beg{thm}\label{thm-uni}
Assume {\bf(H)}, {\bf(A1)}, and there is $\si_0>0$ such that
\beg{align}\label{non-de}
\si(x,\mu)\si^*(x,\mu)\geq \si_0^2,~x\in\R^d,\mu\in\sP^p.
\end{align} 
Suppose that  there is $\ka >0$ such that for all $\mu\in\sP^p$, the invariant probability measure $\cT_\mu$ of $P_t^\mu$ satisfies 
\beg{align}\label{ine-Ta1}
W_p(\nu,\cT_\mu)\leq \sq{2\ka H(\nu|\cT_\mu)},~\nu\in\sP^p.
\end{align}
If $\de<\de_0$ with
\beg{align*}
\de_0&=\sup_{t>t_0, m>m_0} \ff {\si_0(1-\hat C e^{-\hat\la t})K(m,2)[K(m,p)\vee t^{-\ff 1 {p\vee 2}}]} { \si_0K(m,2)+m\sq{\ka t}[K(m,p)\vee t^{-\ff 1 {p\vee 2}}]},\\
t_0&=\hat\la^{-1}{\log\hat C} ,~m_0=\left(\ff {(p-2)^+} {2(p\vee 2)} +  \ff {K_0} 2 \right)^+, 
\end{align*}
then there is a unique stationary probability measure for \eqref{equ_DD}.
\end{thm}

\subsection{Exponential convergence}
To investigate the exponential convergence in  Wasserstein  distance, we  assume that the weak well-posedness of \eqref{equ_DD} holds. Let $P_t^*\mu=\sL_{X_t}^{\P}$ be the law of weak solution with initial distribution $\mu$. Since  \eqref{ine-bsi} implies the pathwise uniqueness of the following equation
\beg{align*}
\d  X_t=b( X_t,P_t^*\mu)\d t+\si(  X_t,P_t^*\mu)\d B_t,
\end{align*}
this equation has a unique strong solution due to the Yamada-Watanabe principle \cite{IW} and the weak well-posedness of \eqref{equ_DD}. As a consequence, \eqref{equ_DD} is strong well-posed.

We first present a result without the assumption {\bf(H)}.
\beg{thm}\label{thm-exp0}
Let $p\geq 2$, and let $\mu\in\sP^q$ with some $q\geq p$. Assume   {\bf (A1)} and  that $|b(0,\cdot)|$ is locally bounded on $\sP^q$ and    \eqref{non-de} holds for some  $\si_0>0$.  We also assume that \eqref{equ_DD} is weak well-posed for $\mu$ and the mapping $t\ra P_t^*\mu$ is locally bounded in $\sP^q$.  Suppose there is a unique stationary distribution $\bar\mu$  for \eqref{equ_DD}, and  there is $\ka>0$ such that \eqref{ine-Ta1} holds with $\cT_\mu $ replaced by $\bar\mu$. If 
\begin{equation}\label{K0-pka}
K_0< \ff {\si_0^2} {2^{3-\ff 4 p}\ka}-\left(\left(\sq{\ff {p-2} p}-\sq{\ff {\si_0^2} {2^{3-\ff 4 p}\ka}}\right)^+\right)^2, 
\end{equation}
then  there is $\de_1>0$ such that for  \eqref{equ_DD} with $\de<\de_1 $,  there exist $\bar C>0,\bar\la>0$ so that 
\beg{align}\label{Mexp}
W_p(P_t^*\mu,\bar\mu) \leq \bar C e^{-\bar\la t }W_p(\mu,\bar\mu).
\end{align}
In particular, if $p=2$, we have
\beg{align}\label{de2}
\de_1&\geq\sq{(2\hat m-K_0)\hat\be^{-1}(\Ph(2)\we\ff 1 2)},\\
\bar\la&\geq \ff {\de^2\hat\be} 2 \left(u-\ff {(1+u)\log (2u) }{\log\ff{2u^2} {\hat\be+(\hat\be-2)u}}\right).\label{con-rat}
\end{align}
where $u=\ff {2\hat m-K_0} {\de^2\hat\be}$, $\hat m= \ff {\si_0^2} {2 \ka}$, $\hat\be=2(1+\ff {\ka \hat m^2} {\si_0^2(2\hat m-K_0)})$ and
\beg{align*}
\Ph(x)=\inf\{v>0 ~|~v (v\hat\be+\hat\be-2)^{-\ff 1 v}\leq x\},~x>0.
\end{align*}

\end{thm}

Instead of the Talagrand inequality for $\bar\mu$, we can also use the twinned Talagrand inequality of $(P_t^{\bar\mu})^*$ and $P_t^*$ to obtain the exponential convergence:
\beg{description}[align=left, noitemsep]
\item[(Ta)] There exist nonempty $\cC\subset\sP$ and $\ka_t>0$  such that $P_t^*\cC\subset \cC$ and
\beg{align}\label{ine-Ta2}
W_p(\nu,(P_t^{\bar\mu})^*\mu)\leq \sq{2\ka_tH(\nu|(P_t^{\bar\mu})^*\mu)},~\nu\in\sP^p,t>0,\mu\in\cC.
\end{align}
\end{description}
The nonempty $\cC$ can not contain all probability measures in $\sP^p$ usually.   See {\bf(A2)}, Example \ref{exa1} and  Lemma \ref{lem-TCI} for concrete conditions  that ensure {\bf(Ta)}.

\beg{thm}\label{thm-exp}
 Let $\mu\in\sP^q$ with some $q\geq p$. The assumption  of Theorem \ref{thm-exp0} hold except  the Talagrand inequality. Assume that {\bf (H)}  holds and there is a unique stationary distribution $\bar\mu$ for \eqref{equ_DD} satisfying  {\bf(Ta)}.  Let
\beg{align*}
\ga(\de,t,\th) =\hat C\left(\ff {1+\th}{\th}\right)^{1-\ff 1 {2\vee p}} \ga_1(\de,t,\th)^{\ff 1 {2\vee p}}e^{- \hat\la t},~\de,\th,t>0,
\end{align*}
with  
\beg{align*}
\ga_1(\de,t,\th) &= 1 +\ff {\de^{2\vee p}C_1(t)} {(1+\th)^{1-2\vee p}} \int_0^t \exp\left\{\int_s^t  \left(\ff {C_1(r)\de^{2\vee p}} {(1+\th)^{1-2\vee p}}+(2\vee p)\hat\la\right) \d r\right\}\d s \\
C_1(t) &= \left(1 +  \ff {t^{\ff {p\vee 2-2 } {2(p\vee 2)}}\sq{\ka_t}(K_0+\ff {(p-2)^+} {p\vee 2}\vee |K_0|)} {2\si_0 \sq{|K_0|\vee \ff {(p-2)^+} {p\vee 2}}}  \right)^{p\vee 2}. 
\end{align*}
Let 
\beg{align*}
\de_2 = \inf\left\{\de >0~\Big|\inf_{t,\th>0}\ga(\de,t,\th)\geq 1\right\}.
\end{align*}
Then $\de_2>0$, and for  \eqref{equ_DD} with $\de<\de_2\we\de_0$,  there are $\bar C>0,\bar\la>0$ such that \eqref{Mexp} holds for every $\mu\in\sP^q\cap\cC$. Let $t_1>0,\th_1>0$ such that $ \ga(\de,t_1,\th_1)<1$. Then $\bar\la\geq   t_1^{-1}\log \ff 1 {\ga(\de,t_1,\th_1)}.$
\end{thm}
\beg{rem}
Recently, exponential convergence in the total variation distance for (reflecting) McKean-Vlasov SDEs has been investigated in \cite[Theorem 2.4]{Wan21b}. The following condition are used  to character the dependence of $b(x,\mu)$ on $\mu$
\beg{align*}
|b(x,\mu)-b(x,\nu)|\leq c\|\mu-\nu\|_{var},~\mu,\nu\in\sP, x\in\R^d,
\end{align*}
where $\|\cdot\|_{\var}$ denotes the  total variation norm. When the constant $c$ is small enough, the existence and uniqueness and the exponential convergence can be established. We adapted a similar argument as in \cite{Wan21b} but investigate the exponential convergence in Wasserstein distance  for equations with distribution dependent $\si$. Since $\si$ here can  be distribution dependent, the coupling used in \cite{Wan21b} can not be applied. 
\end{rem}
\beg{rem}\label{rem0}
If {\bf(A1)} holds with $K_0<0$ and $p\leq 2$, then $m_0=0$ and \eqref{exp-c1} holds with $\hat C=1,~\hat\la=-K_0$. If $\si$ is bounded in addition, then \eqref{ine-Ta1}  holds  with  $\ka=\ff {\|\si\|_\infty} {K_0^-}$, see e.g. \cite[Theorem 5.6]{DGW}. In this situation,  
\beg{align*}
\de_0&=\sup_{t>t_0} \ff {(1-\hat C e^{-\hat\la t})(K(0,2))^2} {  K(0,2)}=\sq{-K_0}.
\end{align*}
Then  $\de<\de_0$ if and only if  $\de^2+K_0<0$, which is the condition ``$C_2>C_1$"  used in  \cite[Theorem 3.1 (2)]{Wan18} for \eqref{Wdis}. Due to $K_0<0$ and $p\leq 2$, $C_1(t)\equiv 1$ and  
\beg{align*}
\inf_{t,\th>0}\ga(\de,t,\th)^{2} & =\inf_{t,\th >0}\left\{  \ff {(1+\th) \left((1+\th)\de^2 e^{(1+\th)\de^2 t}+2\hat\la e^{-2\hat\la t}\right)} {\th ((1+\th)\de^2+2\hat\la)}\right\}  \\
& = \inf_{\th>0} \left\{ \ff {1+\th} {\th}\left(\left(\ff {2\hat\la\de^{-2}} {1+\th}\right)^{\ff {1+\th-2\hat\la\de^{-2}} {1+\th+ 2\hat\la \de^{-2}}}\we 1\right)\right\}.
\end{align*}
Since $2\hat\la/\de_0^2=2$ and 
\beg{align*}
\inf_{\th>0,\de=\de_0} \left\{ \ff {1+\th} {\th}\left(\left(\ff {2\hat\la\de^{-2}} {1+\th}\right)^{\ff {1+\th-2\hat\la\de^{-2}} {1+\th+ 2\hat\la \de^{-2}}}\we 1\right)\right\}=2>1
\end{align*}
we have that $\de_2<\de_0$. We do not obtain sharp estimate for $\de_2$ in the case that $K_0<0$ and $p\leq 2$. 
\end{rem}

\subsection{Corollaries and examples}
Our results can be applied to  SDEs without uniformly dissipative  drifts  investigated by \cite{LuoW,W20}. We use the following condition which is a modified form of \cite[(2.23)]{W20}. In this case, we can derive exponential convergence rate  $\bar \la$   involving $\hat C$ and $\hat \la$.
\beg{description}[align=left, noitemsep]
\item[(A2)] $b$ are continuous  in the  first variable; $\si$ is Lipschitz in the first variable and satisfies \eqref{non-de} and
\beg{align}\label{sup-si}
\|\si\|_{\infty}:=\sup_{x\in\R^d,\mu\in\sP^p}\|\si(x,\mu)\|_{HS}<\infty,
\end{align}  
there exist constant $r_0>0, K_0\geq0, K_1>0$ and $\de\geq 0$ such that $\mu,\nu\in\sP^1$  
\beg{align*} 
&2\<b(x,\mu)-b(y,\nu),x-y\>+\left(2\ff {\|\si\|_\infty^2 } {\si_0^2}-\1_{[d=1]}\right)\|\si(x,\mu)-\si(y,\nu)\|_{HS}^2\nonumber\\
&\qquad\leq \left((K_0+K_1)\1_{[|x-y|\leq r_0]}-K_1\right)|x-y|^2+\de^2W_1(\mu,\nu)^2.
\end{align*}
\end{description}
Then we have the following corollary.
\beg{cor}\label{cor2}
Assume that  {\bf(A2)} holds. Then there are $\hat\la>0$ and $\hat C\geq 1$ such that for $\de<\de_0$ with
\beg{align*}
\de_0&=\sup_{t>\hat\la^{-1}{\log\hat C}} \ff {\si_0(1-\hat C e^{-\hat\la t}) \sq{K_1}} { \sq{2\si_0\|\si\|_\infty\sq{K_1 t}+K_0\|\si\|_\infty^2 t}},
\end{align*}
there is a unique stationary probability measure for \eqref{equ_DD}. 

Suppose  $ |b(0,\cdot)|$ is locally bounded on $\sP^2$. Let $\mu\in\sP$ satisfy 
\beg{align}\label{in-mu0}
\int_{\R^d\times\R^d}e^{\th |x-y|^2 }\mu(\d x)\mu(\d y)<\infty,~\th<\ff {K_1} {4\|\si\|_{\infty}^2}\we \ff {K_3} {4\|\si\|_{\infty}^2}
\end{align}
for some $ K_3>0$, and let  
$$\al=\left(1+\ff {\|\si\|_\infty} {\si_0}\sq{\ff {K_0} {K_1\we K_3}}\right)^{2},\qquad \de_2=\sq{\hat\la\al^{-1}\Ph(2\hat C^2)}$$
where   
$$\Ph(u)=\inf\left\{v>0~|~v^{\ff {v-1} {v+1}}\leq u\right\},~u>0.$$
Then for  \eqref{equ_DD} with $\de<\de_2\we\de_0$,  there are $\bar C>0,\bar\la>0$ such that
\beg{align*}
W_1(P_t^*\mu,\bar\mu) \leq \bar C e^{-\bar\la t }W_1(\mu,\bar\mu).
\end{align*}
with
\beg{align*}
\bar \la &=\left(\ff {\al\de^2+\hat\la} 2\right) \left(\ff {\log(2\hat C^2)} {\log(\hat\la/\al\de^2)}+\ff {\al\de^2-\hat\la}{\al\de^2+\hat\la}\right).
\end{align*} 
\end{cor}

We present the following   example to illustrate {\bf (A2)}, which covers Curie-Weiss mean-field  model  and the granular media model with double-well confinement potential and quadratic interaction, see e.g. \cite{Daw,LWZ,Tug14a}.
\beg{exa}\label{exa1}
Let  $\si$ satisfy \eqref{non-de}, \eqref{sup-si} and there $C_1,\de>0$ so that
\beg{align*}
\left(2\ff {\|\si\|_\infty^2 } {\si_0^2}-\1_{[d=1]}\right)\|\si(x,\mu)-\si(y,\nu)\|_{HS}^2\leq C_1|x-y|^2+\ff {\de^2} 2 W_1(\mu,\nu)^2.
\end{align*}
Let $b$ be of the following form
$$b(x,\mu)=b_1(x)+\int_{\R^d} b_2(x,z)\mu(\d z),~x\in\R^d,$$
where $b_1(x)=-2\th_1 |x|^2x+\ff {\th_2} 2 x$ with some $\th_1>0,\th_2\geq 0$, $b_2:\R^d\times\R^d\ra\R^d$ is continuous and 
\beg{align*}
C_2:=\sup_{(x,z)\in\R^{2d}}|\pp_1b_2(x,z)|<\infty,\qquad \sup_{(x,z)\in\R^{2d}}|\pp_2b_2(x,z)|\leq \ff {\de} {\sq 2}.
\end{align*}
Then
\beg{align*}
2\<b_1(x)-b_1(y),x-y\>\leq - \th_1   |x-y|^4+\th_2|x-y|^2.
\end{align*}
and
\beg{align*}
&2\<\mu(b_2(x,\cdot))-\nu(b_2(y,\cdot)),x-y\>\\
&\qquad =2\<\mu(b_2(x,\cdot)-b_2(y,\cdot))+\mu(b_2(y,\cdot))-\nu(b_2(y,\cdot)),x-y\>\\
&\qquad \leq 2C_2|x-y|^2+ 2\int_{\R^d\times\R^d}\<b_2(y,z_1)-b_2(y,z_2),x-y\>\pi(\d z_1,\d z_2)\\
&\qquad \leq 2C_2|x-y|^2+\sq 2\de W_1(\mu,\nu)|x-y|\\
&\qquad \leq (2C_2+1) |x-y|^2+\ff {\de^2} 2W_1(\mu,\nu)^2.
\end{align*}
Therefore, {\bf(A2)} holds.
\end{exa}


If $\si(x,\mu)$ is independent of $x$, then we can assume
\beg{description}[align=left, noitemsep]
\item[(A2')] $b$ is  continuously differentiable in the  first variable with $\sup_{\mu\in\sP^2}|b(0,\mu)|<\infty$, and there exist   $r_0>0,K_0\geq 0,K_1>0$ and $\de\geq 0$ such that  
\beg{align*}
&2\<b(x,\mu)-b(y,\nu),x-y\>+\|\si( \mu)-\si( \nu)\|_{HS}^2\nonumber\\
&\leq \left((K_0+K_1)\1_{[|x-y|\leq r_0]}-K_1\right)|x-y|^2+\de^2W_2(\mu,\nu)^2,~x,y\in\R^d,~\mu,\nu\in\sP^2.
\end{align*}
\end{description}
Suppose that $\si$ also satisfies \eqref{non-de} and \eqref{sup-si}. Then \eqref{equ_DD} is well-posed by \cite{Wan18}. It follows from \cite[Theorem 2.1(2)]{W20} that there are $\hat\la>0$ and $\hat C\geq 1$ such that \eqref{exp-c1} holds and there is $\ka>0$ such that \eqref{ine-Ta1} holds with $p=2$. Hence, we have the following corollary.
\beg{cor}\label{cor3}
Assume that \eqref{non-de}, \eqref{sup-si} and {\bf(A2')}. Then there are $\hat\la>0$, $\hat C\geq 1$ and $\ka>0$ such that for $\de<\de_0$ with
\beg{align*}
\de_0&=\sup_{t>\hat\la^{-1}{\log\hat C}} \ff {\si_0(1-\hat C e^{-\hat\la t}) } { \sq{\ka(2\si_0\sq t+K_0\ka t)} },
\end{align*}
there is a unique stationary probability measure for \eqref{equ_DD}. If $K_0<\ff {\si_0^2} {2\ka}$ in addition, we obtain the exponential convergence for all $\mu\in\sP^2$ with $\de_1$ and the convergence rate $\bar\la$ given by \eqref{de2} and \eqref{con-rat}.
\end{cor}

\section{Proofs of Proposition \ref{prp1} and theorems}

\noindent{\bf Proof of Proposition \ref{prp1}~}

We prove that  the mapping $\mu\ra \cT_\mu$ is contractive on $\sP^p$. Let $\nu,\mu\in\sP^p$. It follows from \eqref{exp-c1} that
\beg{align}\label{ine-Wp}
W_p(\cT_\mu,\cT_\nu)&\leq W_p(\cT_\mu,(P^\mu_t)^*\cT_\nu)+W_p((P^\mu_t)^*\cT_\nu,\cT_\nu)\nonumber\\
&\leq \hat C e^{-\hat \la t} W_p(\cT_\mu,\cT_\nu)+W_p((P^\mu_t)^*\cT_\nu,(P^\nu_t)^*\cT_\nu).
\end{align}
Let $X_0^\nu=X_0^\mu$ have the law $\cT_\nu$, and let $X_t^\nu$ be the solution of \eqref{equ_dec} with $\mu$ replaced by $\nu$.  If $p\geq 2$, then by the It\^o formula and the H\"older inequality that
\beg{align}\label{Ito-p}
&\d |X_t^\mu-X_t^\nu|^p \nonumber\\
&\quad\leq \ff p 2|X_t^\mu-X_t^\nu|^{p-2}\Big\{2\<b(X_t^\mu,\mu)-b(X_t^\nu,\nu),X_t^\mu-X_t^\nu\>\nonumber\\
&\qquad+\|\si^*(X_t^\mu,\mu)-\si^*(X_t^\nu,\nu)\|_{HS}^2+(p-2)\|\si^*(X_t^\mu,\mu)-\si^*(X_t^\nu,\nu)\|^2\Big\}\d t\nonumber\\
&\qquad +p|X_t^\mu-X_t^\nu|^{p-2}\<X_t^\mu-X_t^\nu,(\si(X_t^\mu,\mu)-\si(X_t^\nu,\nu))\d B_t\>\nonumber\\
&\quad\leq \left(\ff {pK_0} 2|X_t^\mu-X_t^\nu|^{p}+\ff {p\de^2} 2 |X_t^\mu-X_t^\nu|^{p-2}W_p(\mu,\nu)^2\right)\d t\nonumber\\
&\qquad +p|X_t^\mu-X_t^\nu|^{p-2}\<X_t^\mu-X_t^\nu,(\si(X_t^\mu,\mu)-\si(X_t^\nu,\nu))\d B_t\>\nonumber\\
&\quad\leq \left(\ff {pK_0+p-2} 2|X_t^\mu-X_t^\nu|^{p}+ \de^p W_p(\mu,\nu)^p\right)\d t\nonumber\\
&\qquad +p|X_t^\mu-X_t^\nu|^{p-2}\<X_t^\mu-X_t^\nu,(\si(X_t^\mu,\mu)-\si(X_t^\nu,\nu))\d B_t\>.
\end{align}
This implies that
\beg{align*}
W_p((P_t^\mu)^*\cT_\nu,(P_t^\nu)^*\cT_\nu)^p&\leq \E|X_t^\mu-X_t^\nu|^p\\
&\leq \int_0^t e^{-\ff {p(K_0+1)-2} 2 (t-s)} \de^p W_p(\mu,\nu)^p\d s\\
&=\ff {2\de^p t (1-e^{-\ff {p(K_0+1)-2} 2 t})} {p(K_0+1)-2}W_p(\mu,\nu)^p.
\end{align*}
Substituting this into \eqref{ine-Wp}, we arrive at 
\beg{align*}
W_p(\cT_\mu,\cT_{\nu})\leq  \hat C e^{-\hat \la t} W_p(\cT_\mu,\cT_\nu)+\de \left(\ff {2t (1-e^{-\ff {p(K_0+1)-2} 2 t})} {p(K_0+1)-2}\right)^{\ff 1 p}W_p(\mu,\nu).
\end{align*}
Consequently, for $t>\hat \la^{-1}\log \hat C$, 
\beg{align}\label{Wp-de0}
W_p(\cT_\mu,\cT_\nu)&\leq \de \left(\ff {2t (1-e^{-\ff {pK_0+(p-2)} 2 t})} {pK_0+(p-2)}\right)^{\ff 1 p}   (1-\hat C e^{-\hat \la t})^{-1}W_p(\mu,\nu)\nonumber\\
&=\de \left(\ff {2t (1-e^{-(\ff {p\vee 2} 2 {K_0}+\ff {(p-2)^+} 2)t})} {(p\vee 2)K_0+(p-2)^+}\right)^{\ff 1 {p\vee 2}}   (1-\hat C e^{-\hat \la t})^{-1}W_p(\mu,\nu).
\end{align}
If $p\in [1,2)$, then it follows from \eqref{Ito-p} with $p=2$ that
\beg{align*}
\E|X_t^\mu-X_t^\nu|^2&\leq \int_0^t e^{-K_0 (t-s)} \de^2 W_p(\mu,\nu)^2\d s=\ff { \de^2 t (1-e^{ - K_0 t })} {K_0}W_p(\mu,\nu)^2.
\end{align*}
This together with \eqref{ine-Wp} and $W_p\leq W_2$ implies \eqref{Wp-de0}  for $p<2$.  

Therefore, by the definition of $\de_0$,
\beg{align}\label{ine-Wp1}
W_p(\cT_\mu,\cT_\nu)\leq \ff {\de} {\de_0}W_p(\mu,\nu),~\mu,\nu\in\sP^p.
\end{align}
The assertion follows by applying the Banach fixed point theorem to $\cT$ on $\sP^p$.

\qed

In  \eqref{ine-Wp}, we use the synchronous couplings for  $((P^\mu_t)^*\cT_\nu,(P^\nu_t)^*\cT_\nu)$ to estimate $W_p((P^\mu_t)^*\cT_\nu,\cT_\nu)$ under {\bf(A1)}.  We can use the coupling by change of measure and the Talagrand inequality to estimate $W_p((P^\mu_t)^*\cT_\nu,\cT_\nu)$. 

\medskip

\noindent{\bf Proof of Theorem \ref{thm-uni}~}

By \eqref{ine-Wp}, we focus on the estimate of $W_p((P^\mu_t)^*\cT_\nu,\cT_\nu)$. To this aim, we first construct a coupling $(\tld X^\mu_t,X^\nu_t)$ as follows 
\beg{align*}
\d \tld X^{\mu}_t&=b(\tld X^{\mu}_t,\mu)\d t+\si(\tld X^{\mu}_t,\mu)\d B_t\\
&\qquad -m  \si(\tld X^{\mu}_t,\mu)\si^{-1}( X^{\nu}_t,\nu)(\tld X^{\mu}_t-X^{\nu}_t) \d t,~\tld X^\mu_0=\cT_\nu,\\
\d X^{\nu}_t&=b(X^{\nu}_t,\nu)\d t+\si(X^{\nu}_t,\nu)\d B_t,~X^\nu_0=\tld X^\mu_0,
\end{align*}
where $m>m_0$ is a constant and $\tld X^\mu_0=\cT_\nu$  means that $\tld X^\mu_0$  is a random variable with law $\cT_\nu$. By \eqref{ine-bsi} and that $\si$ is bounded and Lipschitz in the first variable, this coupling is well-posed for any $t>0$. Let 
\beg{align*}
\tld B_t=B_t-m\int_0^t \si^{-1}(X^\nu_s,\nu)(\tld X^{\mu}_s-X^{\nu}_s)  \d s.
\end{align*}
Then $(\tld X^\mu_t,X^\nu_t)$ satisfies
\beg{align*}
\d \tld X^{\mu}_t&=b(\tld X^{\mu}_t,\mu)\d t+\si(\tld X^{\mu}_t,\mu)\d \tld B_t,~\tld X^\mu_0=\cT_\nu,\\
\d X^{\nu}_t&=b(X^{\nu}_t,\nu)\d t+\si(X^{\nu}_t,\nu)\d \tld B_t+m(\tld X_t^\mu-X_t^\nu)\d t,~X^\nu_0=\cT_\nu.
\end{align*}
Let 
\beg{align}\label{R_t-mn}
R^{\mu,\nu}_t  =\exp\Big\{&\int_0^t\left\<m\si^{-1}(X^{\nu}_s,\nu)(\tld X^{\mu}_s-X
^{\nu}_s),\d B_s\right\>\nonumber\\
&\qquad -\ff {m^2  } 2\int_0^t \left|\si^{-1}(X^{\nu}_s,\nu)(\tld X^{\mu}_s-X^{\nu}_s)\right|^2 \d s  \Big\}.
\end{align}
As in e.g. \cite{W11}, by combining the Girsanov theorem with the stopping time technique, the Fatou lemma and the martingale convergence theorem, we can derive  that  for any $t>0$, $\{\tld B_s\}_{0\leq s\leq t}$ is a Brownian motion under $\Q=R_t^{\mu,\nu}\P$ and
\beg{align}\label{RlogR0}
\sup_{s\in[0,t]}\E R^{\mu,\nu}_s\log R^{\mu,\nu}_s\leq \ff {t m^2\de^{2} W_p(\mu,\nu)^{2}} {2 \si_0^{2}(2m-K_0)},~t>0.
\end{align} 
It follows from the It\^o formula that 
\beg{align}\label{Ito-mn-2}
\d |\tld X_t^\mu-X_t^\nu|^2&\leq  \left( (K_0-2m)  |\tld X_t^\mu-X_t^\nu|^{2}+ \de^2 W_p(\mu,\nu)^2\right)\d t\nonumber\\
&\qquad +2\<\tld X_t^\mu-X_t^\nu,(\si(\tld X_t^\mu,\mu)-\si(X_t^\nu,\nu))\d \tld B_t\>.
\end{align}
If $p\geq 2$, then it follows from the It\^o formula  that 
\beg{align*}
\d |\tld X_t^\mu-X_t^\nu|^p&\leq  \left(\ff {p(K_0-2m)+p-2} 2|\tld X_t^\mu-X_t^\nu|^{p}+ \de^p W_p(\mu,\nu)^p\right)\d t\nonumber\\
&\qquad +p|\tld X_t^\mu-X_t^\nu|^{p-2}\<\tld X_t^\mu-X_t^\nu,(\si(\tld X_t^\mu,\mu)-\si(X_t^\nu,\nu))\d \tld B_t\>.
\end{align*} 
Combining this with \eqref{Ito-mn-2}, we have that
\beg{align}\label{ine-X-Xp}
\E^\Q |\tld X_t^\mu-X_t^\nu|^{2\vee p}&\leq \int_0^t e^{\ff {(p\vee 2)(K_0-2m)+(p-2)^+} 2 (t-s)}\de^{p\vee 2}W_p(\mu,\nu)^{p\vee 2}\d s\nonumber\\
&\leq \ff {2\de^{p\vee 2}W_p(\mu,\nu)^{p\vee 2}(1-e^{-\ff {(p\vee 2)(2m-K_0)-(p-2)^+} 2})}  {(p\vee 2)(2m-K_0)-(p-2)^+}\nonumber\\
&=\left(\ff {\de W_p(\mu,\nu)}  {K(m,p)\vee t^{-\ff 1 {p\vee 2}}}\right)^{p\vee 2}. 
\end{align}
Denote by $\sL_{X_t^\nu}^\Q$ and $\sL_{\tld X_t^\mu}^\Q$ the law of $X_t^\nu$ and $\tld X_t^\mu$  under $\Q$. Due to the uniqueness of \eqref{equ_dec}, $\sL_{\tld X_t^\mu}^{\Q}=(P^\mu_t)^*\cT_\nu$. Then 
\beg{align*}
\sL_{X_t^\nu}^\Q(f)&=\E^\Q f(X_t^\nu)=\E R_t^{\mu,\nu} f(X_t^\nu)\\
&=\E\left(\E\left[R^{\mu,\nu}_t|X_t^\nu\right] f(X_t^\nu)\right)\\
&=\int_{\R^d}\left(\E\left[R^{\mu,\nu}_t|X_t^\nu=x\right]f(x)\right) \sL_{X_t^\nu}^\P(\d x ),~f\in\sB_b(\R^d).
\end{align*}
Consequently, $\sL_{X_t^\nu}^\Q\ll\sL_{X_t^\nu}^\P$ with
$$\ff {\d\sL_{X_t^\nu}^\Q} {\d\sL_{X_t^\nu}^\P}=\E\left[R^{\mu,\nu}_t|X_t^\nu=x\right],~\sL_{X_t^\nu}^\P\mbox{-a.s.}$$
Hence, it follows from the Jensen inequality and \eqref{RlogR0} that 
\beg{align}\label{ent}
H(\sL_{X_t^\nu}^\Q|\sL_{X_t^\nu}^\P)&=\E\left(\E[R^{\mu,\nu}_t|X_t^\nu]\log\E[R^{\mu,\nu}_t|X_t^\nu]\right)\nonumber\\
&\leq \E\left(\E[R^{\mu,\nu}_t\log R^{\mu,\nu}_t |X_t^\nu]\right)=\E R^{\mu,\nu}_t\log R^{\mu,\nu}_t\nonumber\\
&\leq \ff {m^2 t\de^{2} W_p(\mu,\nu)^{2}} { 2\si_0^{2}(2m-K_0)}=t\left(\ff {m  \de  W_p(\mu,\nu) } { \sq 2\si_0K(m,2)}\right)^2.
\end{align}

By \eqref{ine-X-Xp}, we have that 
\beg{align}\label{ine-Wp-mn0}
W_p((P^\mu_t)^*\cT_\nu,\cT_\nu)&\leq W_p((P^\mu_t)^*\cT_\nu,\sL_{X_t^\nu}^\Q)+W_p(\sL_{X_t^\nu}^\Q,\cT_\nu)\nonumber\\
&\leq \left(\E^\Q|X_t^\mu-X_t^\nu|^{p\vee 2}\right)^{\ff 1 {p\vee 2}}+W_p(\sL_{X_t^\nu}^\Q,\cT_\nu)\nonumber\\
&\leq \ff {\de W_p(\mu,\nu) }  {K(m,p)\vee t^{-\ff 1 {p\vee 2}}}+W_p(\sL_{X_t^\nu}^\Q,\cT_\nu).
\end{align}
It follows from \eqref{ine-Ta1} and \eqref{ent} that 
\beg{align*}
W_p(\sL_{X_t^\nu}^\Q,\cT_\nu)&\leq \sq{2\ka H(\sL_{X_t^\nu}^\Q|\cT_\nu)}=\sq{2\ka H(\sL_{X_t^\nu}^\Q|(P_t^\nu)^*\cT_\nu)}\\
&=\sq{2\ka H(\sL_{X_t^\nu}^\Q|\sL_{X_t^\nu}^\P)}\leq   \ff {m\sq{\ka t}} { \si_0K(m,2)}\de  W_p(\mu,\nu).
\end{align*}
Substituting this into \eqref{ine-Wp-mn0}, we arrive at
\beg{align}\label{ine-W-CH}
W_p((P^\mu_t)^*\cT_\nu,\cT_\nu)&\leq  \left(\ff {1}  {K(m,p)\vee t^{-\ff 1 {p\vee 2}}}+ \ff {m\sq{\ka t}} { \si_0K(m,2)}\right)\de W_p(\mu,\nu).
\end{align}
Taking into account \eqref{exp-c1} and \eqref{ine-Wp}, we have that  for $t>\hat\la^{-1}\log \hat C,~m>m_0$  
\beg{align*}
W_p(\cT_\mu,\cT_\nu) \leq \left(1- \hat C e^{-\hat \la t}\right)^{-1}\left(\ff {1}  {K(m,p)\vee t^{-\ff 1 {p\vee 2}}}+ \ff {m\sq{\ka t}} { \si_0K(m,2)}\right)\de W_p( \mu, \nu).
\end{align*}
Hence, it follows from the definition of $\de_0$ that \eqref{ine-Wp1} holds.  

Therefore, the assertion follows from  \eqref{ine-Wp1} and the Banach fixed point theorem.

\qed

\noindent{\bf Proof of Theorem \ref{thm-exp0}~}

Consider the coupling as follows
\beg{align}\label{tldXt}
\d \hat X_t&=b(\hat X_t,\mu_t)\d t+\si(\hat X_t,\mu_t)\d B_t\nonumber\\
&\qquad -m  \si(\hat X_t,\mu_t)\si^{-1}(\hat X_t^{\bar\mu},\bar\mu)(\hat X_t-\hat X_t^{\bar\mu})  \d t,~\hat X_0=\mu,\\
\d  \hat X_t^{\bar\mu}&=b( \hat X_t^{\bar\mu},\bar\mu)\d t+\si( \hat X_t^{\bar\mu},\bar\mu)\d B_t,~  \hat X_0^{\bar \mu}=\bar\mu,\nonumber
\end{align}
where $m>m_0$ and 
$$\E|\hat  X_0-\hat X_0^{\bar \mu}|^{p}=W_{p}(\mu,\bar\mu)^{p}.$$ 
Since $t\ra P_t^*\mu$ is locally bounded on $\sP^q$, {\bf(A1)} and that $|b(0,\cdot)|$ is locally bounded on $\sP^q$, \eqref{tldXt} is well-posed for any $t>0$. Then $\sL^\P_{\hat X_t^{\bar\mu}}\equiv\bar\mu$ due to {\bf (A1)} which yields the well-posedness of \eqref{equ_dec}.  Let 
\beg{align*}
\hat B_t=B_t-\int_0^t m \si^{-1}( \hat X^{\bar\mu}_s,\bar\mu)(\hat X_s-\hat X_s^{\bar\mu})   \d s.
\end{align*}
Then $(\hat X_t, \hat X_t^{\bar\mu})$ satisfies
\beg{align*}
\d \hat X_t&=b(\hat  X_t,\mu_t)\d t+\si(\hat  X_t,\mu_t)\d \hat B_t,\\
\d \hat X_t^{\bar\mu}&=b( \hat X_t^{\bar\mu},\bar\mu )\d t+\si(\hat X_t^{\bar\mu},\bar\mu )\d \hat B_t+ m (\hat  X_t- \hat X_t^{\bar\mu}) \d t.
\end{align*}
Similarly, as in \cite{W11},  we have that $\{\hat B_s\}_{0\leq s\leq t}$ is a Brownian motion under $\hat\Q=\hat R_t\P$ with
\beg{align*}
\hat R_t  :=\exp\Big\{&\int_0^t\left\<\ m \si^{-1}( \hat X_s^{\bar\mu},\bar\mu)(\hat X_s-\hat X_s^{\bar\mu}),\d B_s\right\>\nonumber\\
&\qquad -\ff {m^2 } 2\int_0^t \left|\si^{-1}(\hat X_s^{\bar\mu},\bar\mu )(\hat X_s- \hat X_s^{\bar\mu})\right|^2 \d s  \Big\}.
\end{align*}
Moreover, by $(p\vee 2)(K_0-2m)+(p-2)^+<0$ since $m > m_0$, we have as \eqref{ine-X-Xp} and \eqref{ent} that  
\beg{align}\label{ine-XXpa}
W_p(\sL^{\hat\Q}_{\hat X_t},\sL^{\hat\Q}_{\hat X^{\bar\mu}_t})^p&\leq \E^{\hat\Q} |\hat X_t-\hat X_t^{\bar\mu}|^{  p}\nonumber\\
&\leq e^{-K(m,p)^{ p} t}\E|\hat X_0-\hat X_0^{\bar \mu}|^{p}+\int_0^t e^{-K(m,p)^{p } (t-s)}\de^{p }W_p(\mu_s,\bar\mu)^{p }\d s\nonumber\\
&\leq e^{-K(m,p)^{p } t}W_{p}(\mu,\bar\mu)^{p } +\de^{p }\int_0^t  W_p(\mu_s,\bar\mu)^{p }\d s
\end{align}
and
\beg{align}\label{h-RlogR}
H(\sL^{\hat\Q}_{\hat X^{\bar\mu}_t}\big|\sL^{\P}_{\hat X^{\bar\mu}_t})&\leq  \E \hat R_t\log \hat R_t=\ff {m^2} 2\E^{\hat \Q}\int_0^t \left|\si^{-1}(\hat X_s^{\bar\mu},\bar\mu )(\hat X_s- \hat X_s^{\bar\mu})\right|^2 \d s\nonumber\\
&\leq \ff {m^2} {2\si_0^2} \int_0^t e^{-(2m-K_0)s} \E|\hat X_0-\hat X_0^{\bar \mu}|^{ 2} \d s\nonumber\\
&\quad +\ff {m^2} {2\si_0^2} \int_0^t \int_0^r e^{-(2m-K_0)(r-s)}\de^2 W_p(\mu_s,\bar \mu)^2\d s\d r\nonumber\\
&\leq \ff {m^2  \left(1-e^{-(2m-K_0)t}\right) } { 2\si_0^{2}(2m-K_0)} W_p(\mu,\bar\mu)^2\nonumber\\
&\quad  +\ff {m^2 \de^{2} } { 2\si_0^{2}(2m-K_0)}\int_0^t W_p(\mu_s,\bar\mu)^{2}\d s.
\end{align}
Combining \eqref{ine-XXpa} with \eqref{h-RlogR}, \eqref{ine-Ta1} for $\bar \mu$ and $\sL^{\P}_{ \hat X_t^{\bar\mu}}\equiv \bar\mu$, we arrive at
\beg{align*}
W_p(P_t^*\mu,\bar\mu)^p&\leq \left(W_p(P_t^*\mu,\sL^{\hat\Q}_{ \hat X_t^{\bar\mu}})+W_p(\sL^{\hat\Q}_{ \hat X_t^{\bar\mu}},\bar\mu)\right)^p\\
&\leq 2^{p-1} W_p(P_t^*\mu,\sL^{\hat\Q}_{ \hat X_t^{\bar\mu}})^p+2^{p-1}W_p(\sL^{\hat\Q}_{ \hat X_t^{\bar\mu}},\bar\mu)^p\\
&\leq 2^{p-1} W_p(\sL^{\hat\Q}_{\hat X_t},\sL^{\hat\Q}_{ \hat X_t^{\bar\mu}})^p+2^{p-1}\left(2\ka H(\sL^{\hat\Q}_{ \hat X_t^{\bar\mu}},\sL^{\P}_{ \hat X_t^{\bar\mu}})\right)^{\ff p 2}\\
&\leq \left(2^{ p-1} e^{-K(m,p)^{p } t}+\ff {2^{\ff {3p} 2-2} \ka^{\ff p 2}  m^p   } { \si_0^{p}(2m-K_0)^{\ff p 2}}\right) W_{p}(\mu,\bar\mu)^{p }\\
&\quad  +2^{p-1}\de^{p }\left(1+\ff { 2^{\ff p 2-1}\ka^{\ff p 2} m^p  t^{\ff p 2-1}} {\si_0^{p}(2m-K_0)^{\ff p 2}}\right)\int_0^t W_p(\mu_s,\bar\mu)^{p}\d s\\
&\equiv a_1(m,t)W_{p}(\mu,\bar\mu)^{p }+\de^p a_2(m,t)  \int_0^t W_p(\mu_s,\bar\mu)^{p}\d s.
\end{align*}
By the Gronwall inequality, we obtain that 
\beg{align*}
W_p(P_t^*\mu,\bar\mu)^p&\leq \ga(\de,m,t)^pW_{p}(\mu,\bar\mu)^{p }
\end{align*}
with
\beg{align*}
\ga(\de,m,t)=\left( a_1(m,t)+\de^p a_2(m,t)  \int_0^t e^{\de^p\int_s^t a_2(m,r)\d r} a_1(m,s)\d s\right)^{\ff 1 p}.
\end{align*}

We choose
$$\hat m=\sq{\ff {\si_0^2} {2^{3-\ff 4 p}\ka}\vee \ff {p-2} p} \sq{\ff {\si_0^2} {2^{3-\ff 4 p}\ka}}.$$
Then by \eqref{K0-pka}, we have that 
\beg{align*}
\left(\ff {\si_0^2} {2^{3-\ff 4 p}\ka}-\ff {\si_0^2} {2^{3-\ff 4 p}\ka}\sq{1-\ff {2^{3-\ff 4 p} \ka} {\si_0^2}K_0}\right)\vee\left(\ff {K_0} 2+\ff {p-2} {2p}\right)\vee 0\\
<\hat m<\ff {\si_0^2} {2^{3-\ff 4 p}\ka} +\ff {\si_0^2} {2^{3-\ff 4 p}\ka}\sq{1-\ff {2^{3-\ff 4 p} \ka} {\si_0^2}K_0}.
\end{align*}
This implies that 
\beg{align*}
\ff {2^{\ff {3p} 2-2}\ka^{\ff p 2} \hat m^p} {\si_0^p(2\hat m-K_0)^{\ff p 2}}<1.
\end{align*}
Then
$$\lim_{t\ra +\infty} a_1(\hat m,t)=\ff {2^{\ff {3p} 2-2}\ka^{\ff p 2} \hat m^p} {\si_0^p(2\hat m-K_0)^{\ff p 2}}<1.$$
Choosing large enough $\hat t$ and small  $\hat \de>0$, one can see that for all $0<\de<\hat \de$, $\ga(\de ,\hat m,\hat t)<1$. Hence 
\beg{align}\label{de2'}
\de_1:=\inf\left\{\de >0~\Big|\inf_{t>0,m>m_0}\ga(\de,m,t)\geq 1\right\}\geq \hat \de>0.
\end{align}
For $\de<\de_1$, there exist $\hat m>0$ and $\hat t>0$ so that $\ga(\de,\hat m,\hat t)<1$ and 
\beg{align*}
W_p(\mu_{\hat t},\bar\mu)\leq \ga(\de,\hat m,\hat t)  W_p(\mu,\bar\mu).
\end{align*}
It follows from the Markov property $P_{t+s}^*=P_t^*P_s^*$ that 
\beg{align*}
W_p(P_t^*\mu,\bar\mu)\leq \left(\ga(\de,\hat m,\hat t)\right)^{\lfloor \ff t { \hat t} \rfloor}W_p(P^*_{t- \hat t\lfloor \ff t { \hat t } \rfloor}\mu,\bar\mu)\leq \bar C(\de,\hat m,\hat t)e^{-\bar\la t }W_p(\mu,\bar\mu),
\end{align*}
where 
\beg{align*}
\bar \la = {\hat t}^{-1}\log \ff 1 {\ga(\de,\hat m,\hat t)},\qquad 
\bar C(\de,\hat m,\hat t)&=\ga(\de,\hat m,\hat t)^{-1}\sup_{0\leq t\leq  \hat t}\ga(\de,\hat m,\hat t).
\end{align*}

In particular, if $p=2$, we have that  $\hat m=\ff {\si_0^2} {2\ka}$ and
\beg{align*}
\ga(\de,\hat m,t)&=2\ff {(2\hat m-K_0)e^{-(2\hat m-K_0)t}} {\de^2\hat\be+2\hat m-K_0}\\
&\quad +2\left(\ff {\de^2\hat\be} {\de^2\hat\be+2\hat m-K_0} +\ff {\ka \hat m^2}{ \si_0^2(2\hat m-K_0)}\right)e^{\de^2\hat\be t},
\end{align*}
Then
\beg{equation*}
\inf_{t>0}\ga(\de,\hat m,t)^2=\beg{cases}
\hat\be,&\de\geq \sq{\ff 1 2(2\hat m-K_0)\hat\be^{-1}},\\
2\left(\ff{\hat\be} 2 + \left(\ff {\hat\be} 2-1\right)u\right)^{\ff u {u+1}} u^{\ff {1-u} {1+u}}, &\de< \sq{\ff 1 2(2\hat m-K_0)\hat\be^{-1}},
\end{cases}
\end{equation*}
where $u=\ff {2\hat m-K_0} {\de^2\hat\be}$ and for $\de< \sq{\ff 1 2(2\hat m-K_0)\hat\be^{-1}}$, the optimal $t$ is 
$$\hat t=\ff {1} {\de^2\be+(2\hat m-K_0)} \log \ff {u^2} {\ff{\hat\be} 2+(\ff{\hat\be} 2-1)u}.$$ 
Thus
\beg{align*}
&\left\{\de>0~\Big|~\inf_{t>0}\ga(\de,\hat m,t)^2\geq 1\right\}\\
&\qquad =\left\{0<\de<\sq{\ff 1 2(2\hat m-K_0)\hat\be^{-1}}~\Big|~2\left(\ff{\hat\be} 2 + \left(\ff {\hat\be} 2-1\right)u\right)^{\ff u {u+1}} u^{\ff {1-u} {1+u}}\geq 1\right\}\\
&\qquad =\left\{0<\de<\sq{\ff 1 2(2\hat m-K_0)\hat\be^{-1}}~\Big|~2\geq  v (v\hat\be+\hat\be-2)^{-\ff 1 v},~v=\ff {\de^2\hat\be} {2\hat m-K_0}\right\}.
\end{align*}
Hence
\beg{align*}
\de_1&\geq \inf\left\{\de>0~\Big|~\inf_{t>0}\ga(\de,\hat m,t)^2\geq 1\right\}\\
& =\inf\left\{\sq{(2\hat m-K_0)\hat \be^{-1}v}~\Big|~0<v<\ff 1 2,~v (v\hat\be+\hat\be-2)^{-\ff 1 v}\leq 2\right\}\\
&= \sq{(2\hat m-K_0)\hat \be^{-1}(\Ph(2)\we \ff 1 2)}.
\end{align*}
Therefore, \eqref{de2} follows, and
\beg{align*}
\bar\la&\geq-\ff 1 {2\hat t}\log\left( 2\left(\ff{\hat\be} 2 + \left(\ff {\hat\be} 2-1\right)u\right)^{\ff u {u+1}} u^{\ff {1-u} {1+u}}\right)\\
&=\ff {\de^2\hat\be} 2 \left(u-\ff {(1+u)\log (2u) }{\log\ff{2u^2} {\hat\be+(\hat\be-2)u}}\right).
\end{align*}

\qed

\noindent{\bf Proof of Theorem \ref{thm-exp}~}

Fix $\mu\in\sP^q\cap\cC$. Since $\bar \mu$ is the stationary distribution,  $\cT_{\bar\mu}=\bar\mu$. Since \eqref{equ_DD} is weak well-posed for $\mu\in\sP^q$, $P_t^*\mu$ is well-defined.  By \eqref{exp-c1}, we have that 
\beg{align}\label{ine-Wmu0}
W_p(P_t^*\mu, \bar\mu)&\leq W_p(P_t^*\mu,(P_t^{\bar\mu})^*\mu)+W_p((P_t^{\bar\mu})^*\mu,\bar\mu)\nonumber\\
&\leq W_p(P_t^*\mu,(P_t^{\bar\mu})^*\mu)+\hat C  e^{-\hat\la t} W_p(\mu,\bar\mu).
\end{align}
Denote $\mu_t=P_t^*\mu$. We consider the following coupling
\beg{align} 
\d \tld X_t&=b(\tld X_t,\mu_t)\d t+\si(\tld X_t,\mu_t)\d B_t\nonumber\\
&\qquad -m  \si(\tld X_t,\mu_t)\si^{-1}(X_t^{\bar\mu},\bar\mu)(\tld X_t-X_t^{\bar\mu})  \d t,~\tld X_0=\mu,\\
\d  X_t^{\bar\mu}&=b( X_t^{\bar\mu},\bar\mu)\d t+\si( X_t^{\bar\mu},\bar\mu)\d B_t,~  X_0^{\bar \mu}=\tld X_0,\nonumber
\end{align}
where $m>m_0$.  Let 
\beg{align*}
\tld B_t=B_t-\int_0^t m \si^{-1}( X^{\bar\mu}_s,\bar\mu)(\tld X_s- X_s^{\bar\mu})   \d s.
\end{align*}
Then $(\tld X_t, X_t^{\bar\mu})$ satisfies
\beg{align*}
\d \tld X_t&=b(\tld X_t,\mu_t)\d t+\si(\tld X_t,\mu_t)\d \tld B_t,\\
\d   X_t^{\bar\mu}&=b( X_t^{\bar\mu},\bar\mu )\d t+\si( X_t^{\bar\mu},\bar\mu )\d \tld B_t+ m (\tld X_t- X_t^{\bar\mu}) \d t.
\end{align*}
We can prove that $\{\tld B_s\}_{0\leq s\leq t}$ is a Brownian motion under $\tld\Q=\tld R_t\P$ with
\beg{align}\label{eq-R_t}
\tld R_t  :=\exp\Big\{&\int_0^t\left\<\ m \si^{-1}( X_s^{\bar\mu},\bar\mu)(\tld X_s-X_s^{\bar\mu}),\d B_s\right\>\nonumber\\
&\qquad -\ff {m^2 } 2\int_0^t \left|\si^{-1}( X_s^{\bar\mu},\bar\mu )(\tld X_s- X_s^{\bar\mu})\right|^2 \d s  \Big\},
\end{align}
and
\beg{align}\label{RlogR1}
\sup_{s\in[0,t]}\E \tld R_s\log \tld R_s&\leq \ff {m^2 \de^{2} } { 2\si_0^{2}(2m-K_0)}\int_0^t W_p(\mu_s,\bar\mu)^{2}\d s,\\
\E^{\tld \Q} |X_t^{\bar\mu}-\tld X_t|^{2\vee p}&\leq \int_0^t e^{ -K(m,p)^{2\vee p} (t-s)}\de^{p\vee 2}W_p(\mu_s,\bar\mu)^{p\vee 2}\d s\nonumber\\
&\leq \de^{p\vee 2}\int_0^t  W_p(\mu_s,\bar\mu)^{p\vee 2}\d s, ~t>0,\nonumber
\end{align}
where we have used in the last inequality that  $(p\vee 2)(K_0-2m)+(p-2)^+\leq 0$. It follows from the uniqueness in law of \eqref{equ_DD} that  $\sL_{\tld X_t}^{\tld\Q}=\mu_t$.   Then
\beg{align*}
W_p(P_t^*\mu,\sL^{\tld\Q}_{ X_t^{\bar\mu}})&\leq W_{p\vee 2}(\sL_{\tld X_t}^{\tld\Q},\sL^{\tld\Q}_{ X_t^{\bar\mu}})\leq \left(\E^{\tld\Q}| X_t^{\bar\mu}-\tld X_t |^{p\vee 2}\right)^{\ff 1 {p\vee 2}}\\
&\leq \de \left(\int_0^t  W_p(\mu_s,\bar\mu)^{p\vee 2}\d s\right)^{\ff 1 {p\vee 2}}.
\end{align*}
Due to $\sL_{ X_t^{\bar\mu}}^{\tld\Q}(f) =\E \tld R_t f( X_t^{\bar\mu})$ and \eqref{RlogR1}, we also have by $\mu\in\cC$ and \eqref{ine-Ta2} that
\beg{align*}
W_p(\sL_{ X_t^{\bar\mu}}^{\tld\Q},\sL^\P_{ X_t^{\bar\mu}})&\leq  \sq{2\ka_t H(\sL_{X_t^{\bar\mu}}^{\tld\Q}|\sL^\P_{ X_t^{\bar\mu}})}\leq \sq{2\ka_t\E \tld R_t\log \tld R_t}\\
&\leq \ff {m \de \sq{\ka_t} } { \si_0  \sq{2m-K_0}}\left(\int_0^t W_p(\mu_s,\bar\mu)^{2}\d s\right)^{\ff 1 2}\\
&\leq \ff {m \de \sq{\ka_t} t^{\ff {p\vee 2-2 } {2(p\vee 2)}} } { \si_0  \sq{2m-K_0}}\left(\int_0^t W_p(\mu_s,\bar\mu)^{2\vee p}\d s\right)^{\ff 1 {2\vee p}}.
\end{align*}
Hence,
\beg{align}\label{ine-WXY}
W_p(P_t^*\mu,(P_t^{\bar\mu})^*\mu)& \leq W_p(P_t^*\mu,\sL^{\tld\Q}_{ X_t^{\bar\mu}})+W_p(\sL^{\tld\Q}_{ X_t^{\bar\mu}},(P_t^{\bar\mu})^*\mu)\nonumber\\
&\leq \left(\E^{\tld\Q} |\tld X_t- X_t^{\bar\mu} |^{2\vee p}\right)^{\ff 1{2\vee p}}+W_p(\sL_{X_t^{\bar\mu}}^{\tld\Q},\sL^\P_{  X_t^{\bar\mu}})\nonumber\\
&\leq \de\left(1 + \ff {m  \sq{\ka_t} t^{\ff {p\vee 2-2 } {2(p\vee 2)}}} { \si_0  \sq{2m-K_0}}\right)\left(\int_0^t  W_p(\mu_s,\bar\mu)^{p\vee 2}\d s\right)^{\ff 1 {p\vee 2}}.
\end{align}
Combining this with \eqref{ine-Wmu0}, it follows from the Jensen inequality that for all $\th>0$
\beg{align*}
W_p(\mu_t,\bar\mu)^{2\vee p}&\leq  \left(W_p(P_t^*\mu,(P_t^{\bar\mu})^*\mu)+\hat C e^{-\la t} W_p(\mu,\bar\mu)\right)^{2\vee p}\\
&=  \left(\ff {1+\th} {1+\th}W_p(P_t^*\mu,(P_t^{\bar\mu})^*\mu)+\ff {\th} {1+\th}\ff {1+\th} {\th}\hat C e^{-\la t} W_p(\mu,\bar\mu)\right)^{2\vee p}\\
&\leq (1+\th)^{2\vee p-1} \de^{2\vee p}\left(1 + \ff {m  \sq{\ka_t} t^{\ff {p\vee 2-2 } {2(p\vee 2)}}} { \si_0  \sq{2m-K_0}}\right)^{p\vee 2} \int_0^t  W_p(\mu_s,\bar\mu)^{p\vee 2}\d s\\
&\quad +\left(\ff {1+\th} {\th}\right)^{2\vee p-1}\hat C^{2\vee p}e^{- (2\vee p )\hat \la t} W_p(\mu,\bar\mu)^{2\vee p},~t>0.
\end{align*}
Let 
\beg{align*}
\ga(\de,t,m,\th) = \hat C \left(\ff {1+\th} {\th}\right)^{1-\ff 1 {2\vee p}}\ga_1(\de,t,m,\th)^{\ff 1 {2\vee p}}e^{- \hat\la t},~\de,\th,t>0,
\end{align*}
with  
\beg{align*}
\ga_1(\de,t,m,\th) &= 1+\ff {\de^{2\vee p}C_1(t,m) } {(1+\th)^{1-2\vee p}}\int_0^t \exp\left\{\int_s^t  \left(\ff {C_1(r,m) \de^{2\vee p}} {(1+\th)^{1-2\vee p}}+(2\vee p)\hat\la\right) \d r\right\}\d s \\
C_1(t,m) &= \left(1 + \ff {m  \sq{\ka_t} t^{\ff {p\vee 2-2 } {2(p\vee 2)}}} { \si_0  \sq{2m-K_0}}\right)^{p\vee 2}. 
\end{align*}
It follows from the Gronwall inequality that
\beg{align*}
W_p(\mu_t,\bar\mu)^{2\vee p}\leq \ga(\de,t,m,\th)^{2\vee p}W_p(\mu,\bar\mu)^{2\vee p}.
\end{align*}
Taking optiaml $m$, we have that
\beg{align}\label{optm}
\inf_{m>m_0}\ff m {\sq{2m-K_0}}&=\ff m {\sq{2m-K_0}}\Big|_{m=\ff 1 2\left(\ff {(p-2)^+} {p\vee 2}\vee |K_0|+K_0\right)}\nonumber\\
&= \ff {K_0+\ff {(p-2)^+} {p\vee 2}\vee |K_0|} {2\sq{|K_0|\vee \ff {(p-2)^+} {p\vee 2}}}.
\end{align}
Then
\beg{align*}
\inf_{t,\th>0,m>m_0}\ga(\de,t,m,\th)&=\inf_{t,\th>0}\ga\left(\de,t,\ff 1 2\left(\ff {(p-2)^+} {p\vee 2}\vee |K_0|+K_0\right),\th\right)\\
&\equiv \inf_{t,\th>0}\ga(\de,t,\th).
\end{align*}
Hence
\beg{align*}
\de_2 &=\inf\left\{\de >0~\Big|\inf_{t,\th>0}\ga(\de,t,\th)\geq 1\right\}\\
& = \inf\left\{\de >0~\Big|\inf_{t,\th>0,m>m_0}\ga(\de,t,m,\th)\geq 1\right\}.
\end{align*}
Since $\displaystyle \lim_{t\ra+\infty}\lim_{\de\ra 0^+}\ga(\de,t,m,\th)=0$, there are $t',m'>m_0,\th'>0$ and $\de'>0$ such that for all $\de<\de'$ it holds that $\ga(\de,t,m,\th)<1$. Thus $\de_2\geq \de'>0$.

For $\de<\de_2\we \de_0$, let $(t_1 , \th_1)$ so that $\ga(\de,t_1,\th_1)<1$. Then
\beg{align*}
W_p(\mu_{t_1},\bar\mu)\leq \ga(\de,t_1,\th_1)  W_p(\mu,\bar\mu).
\end{align*}
Due to {\bf(Ta)} and $P_t^*\mu\in\sP^q$, we have that $P_t^*(\sP^q\cap\cC)\subset \sP^q\cap\cC$. Then, it follows from the Markov property $P_{t+s}^*=P_t^*P_s^*$  that 
\beg{align*}
W_p(P_t^*\mu,\bar\mu)\leq \left(\ga(\de,t_1,\th_1)\right)^{\lfloor \ff t { t_1} \rfloor}W_p(P^*_{t- t_1\lfloor \ff t { t_1} \rfloor}\mu,\bar\mu)\leq \bar C(\de,t_1,\th_1)e^{-\bar\la t }W_p(\mu,\bar\mu),
\end{align*}
where $\bar \la = t_1^{-1}\log \ff 1 {\ga(\de,t_1,\th_1)}$ and
\beg{align*}
\bar C(\de,t_1,\th_1)&=\ga(\de,t_1,\th_1)^{-1}\sup_{0\leq t\leq   t_1}\ga(\de,t,\th_1).
\end{align*}

\qed

\section{Proofs of corollaries}

\noindent{\bf Proof of Corollary \ref{cor2}~}

To prove Corollary \ref{cor2}, we first establish the $W_1$-transportation cost inequality under the assumption {\bf (A2)}. Let $\nu_0\in\sP^1$ and $\{\nu_t\}_{t\geq 0}\subset \sP^1$ such that the following SDE has a unique solution 
\beg{align*}
\d Y_t=b(Y_t,\nu_t)\d t+\si(Y_t,\nu_t)\d B_t,~Y_0=\nu_0.
\end{align*}

\beg{lem}\label{lem-TCI}
Assume that $b,\si$ satisfies {\bf(A2)}. If $\nu_0$ satisfies \eqref{in-mu0}, 
then   $\sL_{Y_t}^{\P}$ satisfies \eqref{in-mu0} and $W_1$-transportation cost inequality 
\beg{align}\label{ine-W1-entr}
W_1(\nu,\sL^{\P}_{Y_t})\leq \sq{2\|\si\|_{\infty}^2( K_1\we K_3)^{-1}H(\nu|\sL^{\P}_{Y_t})},~t\geq 0.
\end{align}
\end{lem}
\beg{proof}
We establish $W_1$-transportation cost inequality for \eqref{equ_dec} by using \cite[(1.5) and Theorem 2.3]{DGW}, see also \cite[Theorem 3.2]{SYZ}. 

Let $\{\bar B_t\}_{t\geq 0}$ be a Brownian motion independent of $\{B_t\}_{t\geq 0}$, and let $\bar Y_t$ be the solution of the following equation
\beg{align*}
\d \bar Y_t=b(\bar Y_t,\nu_t)\d t+\si(\bar Y_t,\nu_t)\d \bar B_t,~\bar Y_0=\mu_0,
\end{align*}
and $\bar Y_0$ is independent of $Y_0$. It follows from the It\^o formula that 
\beg{align*}
\d |Y_t-\bar Y_t|^2&=2\<Y_t-\bar Y_t,b(Y_t,\nu_t)-b(\bar Y_t,\nu_t)\>\d t\\
&\quad +(\|\si(Y_t,\nu_t)\|_{HS}^2+\|\si(\bar y_t,\nu_t)\|_{HS}^2)\d t\\
&\quad +2\<Y_t-\bar Y_t,\si(Y_t,\nu_t)\d B_t-\si(\bar Y_t,\nu_t)\d\bar B_t\>\\
&\leq \left((K_0+K_1)\1_{[|Y_t-\bar Y_t|\leq r_0]}-K_1\right)|Y_t-\bar Y_t|^2\d t\\
&\quad -\|\si(Y_t,\nu_t)-\si(\bar Y_t,\nu_t)\|_{HS}^2\d t\\
&\quad +(\|\si(Y_t,\nu_t)\|_{HS}^2+\|\si(\bar y_t,\nu_t)\|_{HS}^2)\d t\\
&\quad +2\<Y_t-\bar Y_t,\si(Y_t,\nu_t)\d B_t-\si(\bar Y_t,\nu_t)\d\bar B_t\>\\
&\leq \left((K_0+K_1)\1_{[|Y_t-\bar Y_t|\leq r_0]}-K_1\right)|Y_t-\bar Y_t|^2\d t+2\|\si\|_{\infty}^2\d t\\
&\quad +2\<Y_t-\bar Y_t,\si(Y_t,\nu_t)\d B_t-\si(\bar Y_t,\nu_t)\d\bar B_t\>.
\end{align*}
This implies that
\beg{align*}
\d e^{\th|Y_t-\bar Y_t|^2}&\leq \th e^{\th |Y_t-\bar Y_t|^2}\Big\{\left((K_0+K_1)\1_{[|Y_t-\bar Y_t|\leq r_0]}-K_1+4\th\|\si\|_{HS,\infty}^2\right)|Y_t-\bar Y_t|^2\d t\\
&\quad +2\|\si\|_{\infty}^2\d t+2\<Y_t-\bar Y_t,\si(Y_t,\nu_t)\d B_t-\si(\bar Y_t,\nu_t)\d\bar B_t\>\Big\}.
\end{align*}
For every $\th<\ff {K_1} {4\|\si\|_{\infty}^2}$, we have that $K_1-4\th\|\si\|_{\infty}^2>0$. Then there are  positive constants $\tld K_0$ and $\tld K_1$ independent of $t$ such that
\beg{align*}
e^{\th x^2}\left\{\left((K_0+K_1)\1_{[x\leq r_0]}-K_1+4\th\|\si\|_{\infty}^2\right)x^2+2\|\si\|_{\infty}^2\right\}\leq \tld K_0-\tld K_1e^{\th x^2}.
\end{align*}
Thus
\beg{align*}
\d e^{\th|Y_t-\bar Y_t|^2}&\leq \left(\tld K_0-\tld K_1 e^{\th|Y_t-\bar Y_t|^2}\right)\d t\\
&\quad + 2e^{\th|Y_t-\bar Y_t|^2}\<Y_t-\bar Y_t,\si(Y_t,\mu_1)\d B_t-\si(\bar Y_t,\mu_1)\d\bar B_t\>.
\end{align*}
This yields that
\beg{align}\label{in-adY}
\E e^{\th|Y_t-\bar Y_t|^2} &\leq e^{-\th \tld K_1 t}\E e^{\th |Y_0-\bar Y_0|^2}+\tld K_0 \int_0^t e^{-\th \tld K_1(t-s)}\d s\nonumber\\
&\leq e^{-\th \tld K_1 t} \int_{\R^d\times\R^d}e^{\th |x-y|^2}\nu_0(\d x)\nu_0(\d y)+\ff { \tld K_0} {\th \tld K_1},~t\geq 0.
\end{align}
Due to the uniqueness of solutions to \eqref{equ_dec} and that $Y_0$ and $\bar Y_0$ are independent with the same law, $Y_t$ and $\bar Y_t$ are independent with the same law. Then it follows from \eqref{in-adY} that for any $\th<\ff {K_1} {4\|\si\|_{\infty}^2}\we \ff {K_3} {4\|\si\|_{\infty}^2}$
\beg{align*}
\int_{\R^d\times\R^d} e^{\th |y_1-y_2|^2}\sL_{Y_t}^{\P}(\d y_1)\sL_{Y_t}^{\P}(\d y_2)=\E e^{\th|Y_t-\bar Y_t|^2} <+\infty,~t\geq 0.
\end{align*}
Hence, according to \cite[(1.5) and Theorem 2.3]{DGW},  \eqref{ine-W1-entr} holds.

\end{proof}

We now turn to the proof of Corollary \ref{cor2}.  According to  \eqref{non-de}, \eqref{sup-si} and \cite[Theorem 2.6 and (2.23)]{W20},   {\bf (H)} holds with $p=1$. It is clear that {\bf (A2)} implies that {\bf (A1)} holds with $p=1$.

(1) Due to {\bf (H)}, $(P_t^\mu)^*{\de_x}\ra \cT_\mu$ weakly as $t\ra +\infty$, where $\de_x$ is the Dirac measure centred on $x\in\R^d$.  This, together with that  \eqref{in-mu0} holds for $\mu_0=\de_x$ and any $K_3>0$, implies by Lemma \ref{lem-TCI} and \cite[Lemma 2.2]{DGW} that \eqref{ine-Ta1} holds for any $\cT_\mu$ with $\ka=\|\si\|_\infty^2 K_1^{-1}$. Then by  Theorem \ref{thm-uni}
\beg{align*}
&\sup_{t>t_0, m>m_0} \ff {\si_0(1-\hat C e^{-\hat\la t})K(m,2)[K(m,p)\vee t^{-\ff 1 {p\vee 2}}]} { \si_0K(m,2)+m\sq{\ka t}[K(m,p)\vee t^{-\ff 1 {p\vee 2}}]}\\
&\quad \geq \sup_{t>\hat\la^{-1}{\log\hat C}, m>K_0/2} \ff {\si_0(1-\hat C e^{-\hat\la t})\sq{2m-K_0} } { \si_0 +m\ff {\|\si\|_\infty} {\sq{K_1}} \sq t}\\
&\quad =\sup_{t>\hat\la^{-1}{\log\hat C}} \ff {\si_0(1-\hat C e^{-\hat\la t})\sq{K_1 } } { \sq{2\si_0\|\si\|_\infty\sq{K_1 t} +\|\si\|_\infty^2K_0 t}},
\end{align*}
where in the last equality we  set  
$$m=\ff 1 2\left\{\left(\ff {2\si_0\sq{K_1}+\|\si\|_\infty K_0\sq t} {\|\si\|_\infty\sq t}\right)^2+K_0\right\}.$$
Then we have proved the first assertion. 

(2) It follows from Lemma \ref{lem-TCI} and \eqref{in-mu0} that \eqref{ine-Ta2} holds with $\ka_t\equiv \|\si\|_{\infty}^2( K_1\we K_3)^{-1}$ and $\cC$  consists of probability measures satisfying \eqref{in-mu0}. The inequality \eqref{in-mu0} also yields that $\mu\in\sP^2$. Since $W_1\leq W_2$, {\bf(A2)} implies the strong well-posedness of \eqref{equ_DD}  with initial distribution $\mu$ and $t\ra \P_t^*\mu$ is locally bounded in $\sP^2$, see e.g. \cite{Wan18}.  It follows from Theorem \ref{thm-exp} with $\th=1$ and Lemma \ref{lem-TCI} that
\beg{align*}
\ga(\de,t,1)^2=\ff {2\hat C^2} {\al \de^2+\hat\la}\left(\al\de^2 e^{2\al\de^2t}+\hat\la e^{-2\hat\la t}\right)
\end{align*}
Then
\beg{equation*}
\inf_{t>0}\ga(\de,t,1)^2=\beg{cases}
 \ff {2\hat C^2\al\de^2} {\al \de^2+\hat\la}, &\de> \sq{\hat \la \al^{-1}}\\
 2\hat C^2\left(\ff {\hat\la} {\al \de^{2}}\right)^{\ff {\al\de^{2}-\hat\la} {\al\de^2+\hat\la}},&\de\leq \sq{\hat \la \al^{-1}}
\end{cases}.
\end{equation*}
Since for $\de> \sq{\hat \la \al^{-1}}$
$$\ff {2\hat C^2\al\de^2} {\al \de^2+\hat\la}>\hat C^2\geq 1,$$
we have that
\beg{align*}
\inf&\left\{\de >0~\Big|\inf_{t,\th>0}\ga(\de,t,\th)\geq 1\right\}\\
&\geq \inf\left\{\de >0~\Big|\inf_{t >0}\ga(\de,t,1)\geq 1\right\}\\
&=\inf\left\{0<\de\leq \sq{\hat \la \al^{-1}}~\Big|~2\hat C^2\geq \left(\ff {\al \de^{2}}  {\hat\la} \right)^{\ff {\al\de^{2}-\hat\la} {\al\de^2+\hat\la}} \right\}\\
&=\inf\left\{\de~\Big|~0<u\leq 1,~u^{\ff {u-1} {u+1}}\leq 2\hat C^2,~u= \ff {\al \de^{2}} {\hat\la}\right\}\\
&=\sq{\hat\la \al^{-1} \Ph(2\hat C^2)},
\end{align*}
where in the last equality we have used that the function $u^{\ff {u-1} {u+1}}$ decrease in $(0,1)$ and increase in $(1,+\infty)$. Then we obtain $\de_2$.

For $\de<\de_2\we\de_0$, then $\de<\sq{\hat\la \al^{-1}\Ph(2\hat C^2)}$ which yields that $\hat\la/(\al\de^2)>1$. Thus we can choose $\hat t =\ff {\log (\hat\la/\al\de^2)} {\al\de^2+\hat\la}$ so that  
$$\ga(\de,\hat t(\de),1)^2=2\hat C^2\left(\ff {\hat\la} {\al \de^{2}}\right)^{\ff {\al\de^{2}-\hat\la} {\al\de^2+\hat\la}}<1$$
Then 
\beg{align*}
W_1(\mu_{\hat t },\mu_1)&\leq  \ga(\de,\hat t ,1) W_1(\mu,\mu_1).
\end{align*}
It follows from the Markov property $P_{t+s}^*=P_t^*P_s^*$ that 
\beg{align*}
W_1(P_t^*\mu,\mu_1)\leq \left(\ga(\de,\hat t ,1)\right)^{\lfloor \ff t {\hat t } \rfloor}W_1(P^*_{t-\hat t \lfloor \ff t {\hat t } \rfloor},\mu_1)\leq \bar C e^{-\bar\la t }W_1(\mu,\mu_1),
\end{align*}
with
\beg{align*}
\bar \la&=\ff 1 {\hat t }\log \ff 1 {\ga(\de,\hat t ,1)}=\left(\ff {\al\de^2+\hat\la} 2\right) \left(\ff {\log(2\hat C^2)} {\log(\hat\la/\al\de^2)}+\ff {\al\de^2-\hat\la}{\al\de^2+\hat\la}\right),\\
\bar C &=\ga(\de,\hat t ,1)^{-1}\sup_{0\leq t\leq \hat t }\ff {C_1^2\left(\al\hat \de^2 e^{2\al\de^2 t}+\hat \la  e^{-2\hat \la t}\right)} {\al\de^2+\hat \la}.
\end{align*}

\qed

\noindent{\bf Proof of Corollary \ref{cor3}~}

The condition {\bf (A2')} yields the strong well-posedness to \eqref{equ_DD}, see e.g. \cite{Wan18}. It follows from {\bf(A2')} and \cite[Corollary 1.8]{LuoW} that for any $\mu\in\sP^2$, there is $\cT_\mu\in\cap_{q\geq 1}\sP^q$ and \eqref{exp-c1} holds for $p=1$.  We use \cite[Theorem 2.1 (2)]{W20} to prove that  \eqref{exp-c1}  and   \eqref{ine-Ta1} holds  for $p=2$.  To this end, we only need to prove that $P_t^\mu$ has an invariant probability measure $\cT_\mu$ and  there is $\th>0$ independent of $\mu$ such that $\cT_\mu(e^{\th|\cdot|^2})<\infty$.  For the solution to \eqref{equ_dec} with $X^\mu_0=0$, it follows from {\bf (A2')}, the It\^o formula and the H\"older inequality that
\beg{align*}
\d |X_t^\mu|^2&\leq \left((K_0+K_1)\1_{[|X_t^\mu|\leq r_0]}-K_1\right)|X_t^\mu|^2+2\<X_t^\mu,\si(\mu)\d B_t\>\\
&\quad +\|\si(\mu)\|_{HS}^2\d t+2\<b(0,\mu),X_t^\mu\>\d t\\
&\leq \left((K_0+K_1)\1_{[|X_t^\mu|\leq r_0]}-\ff {K_1} 2\right)|X_t^\mu|^2+2\<X_t^\mu,\si(\mu)\d B_t\>\\
&\quad +\|\si(\mu)\|_{HS}^2\d t+  K_1^{-1}\sup_{\mu\in\sP^2}|b(0,\mu)|^2\d t.
\end{align*}
Then
\beg{align*}
\d e^{\th |X_t^\mu|^2}&\leq \th e^{\th |X_t^\mu|^2}  \left((K_0+K_1)\1_{[|X_t^\mu|\leq r_0]}-\ff {K_1} 2\right)|X_t^\mu|^2\d t\\
&\quad +\left\{\|\si\|_\infty^2+\ff {\sup_{\mu\in\sP^2}|b(0,\mu)|^2} {K_1}\right\}\d t\\
&\quad +2\th^2 e^{\th |X_t^\mu|^2}\|\si\|_\infty^2|X_t^\mu|^2\d t+2\th e^{\th |X_t^\mu|^2}\<X_t^\mu,\si(\mu)\d B_t\>.
\end{align*}
This implies that there is $\th_0>0$ independent of $\mu$ such that for  $\th<\th_0$
\beg{align*}
\d e^{\th |X_t^\mu|^2}\leq \left(C_1-C_2   e^{\th |X_t^\mu|^2}\right)\d t+2\th e^{\th |X_t^\mu|^2}\<X_t^\mu,\si(\mu)\d B_t\>
\end{align*}
with some $C_1,C_2>0$ independent of $\mu$. Then 
$$\ff 1 t\int_0^t \E e^{\th |X_s^\mu|^2}\d s\leq \ff {C_1} {C_2},$$
which implies $\sup_{\mu\in \sP^2}\cT_\mu(e^{\th |\cdot|^2})<\infty$. From this and {\bf (A2')}, which implies that
\beg{align*}
2\<b(x,\mu)-b(y,\mu),x-y\>\leq \left((K_0+K_1)\1_{[|x-y|\leq r_0]}-K_1\right)|x-y|^2,~\mu\in\sP^2,
\end{align*}
we have that    \eqref{exp-c1}  and   \eqref{ine-Ta1} holds  for $p=2$ according to  \cite[Theorem 2.1 (2)]{W20}.


Applying Theorem \ref{thm-uni}  with $p=2$, we have that 
\beg{align*}
&\sup_{t>t_0, m>m_0} \ff {\si_0(1-\hat C e^{-\hat\la t})K(m,2)[K(m,p)\vee t^{-\ff 1 {p\vee 2}}]} { \si_0K(m,2)+m\sq{\ka t}[K(m,p)\vee t^{-\ff 1 {p\vee 2}}]}\\
&\quad\geq \sup_{t>t_0, m>m_0} \ff {\si_0(1-\hat C e^{-\hat\la t})\sq{2m-K_0}} { \si_0 +m\sq{\ka t}}\\
&\quad = \sup_{t>t_0} \ff {\si_0(1-\hat C e^{-\hat\la t})} {\sq{\ka(2\si_0\sq t+K_0\ka t)}}.
\end{align*}
Hence the first assertion follows. 

Applying Theorem \ref{thm-uni}  with $p=2$, we prove the second assertion.

\qed




\bigskip
 
\noindent\textbf{Acknowledgements}

\medskip

The  author was supported by  the disciplinary development project of Central University of Finance and Economics, and the National Natural Science Foundation of China (Grant No. 11901604, 11771326).


\end{document}